\documentclass[12pt]{amsart}
\usepackage{amsfonts}
\usepackage{amssymb}
\usepackage{amsmath}
\usepackage{amsfonts}
\usepackage{graphicx}
\usepackage{amscd}
\usepackage{rotating}
\usepackage{amsmath,amsfonts,amssymb,amsthm,mathrsfs,xypic}
\xyoption{curve}
\usepackage{fancybox}
\usepackage{color}
\usepackage[perpage,symbol]{footmisc}
\setfnsymbol{wiley}

\usepackage{amsfonts}
\usepackage{graphicx}
\usepackage{amscd}
\usepackage[all]{xy}
\usepackage{amsmath,accents,amsfonts,amssymb,amsthm,mathrsfs,xypic,extarrows,mathtools}
\usepackage{leftidx}
\newtheorem{proposition}{Proposition}[section]
\newtheorem{lemma}[proposition]{Lemma}
\newtheorem{corollary}[proposition]{Corollary}
\newtheorem{theorem}[proposition]{Theorem}

\theoremstyle{definition}
\newtheorem{definition}[proposition]{Definition}

\theoremstyle{remark}

\newcommand{\Hom}{\mbox{\rm Hom}}

\newcommand{\im}{\mbox{\rm Im}}

\newcommand{\Diag}{\mbox{\rm Diag}}

\newcommand{\Gpd}{\mbox{\rm Gpd}}

\begin{document}

\title{Stability of Gorenstein objects in triangulated categories}
\author{ZHANPING WANG \ \ \ \ CHUNLI LIANG}

\footnote[0]{Address
correspondence to Zhanping Wang, Department of Mathematics, Northwest Normal University, Lanzhou 730070, P.R. China.}\footnote[0]{E-mail: wangzp@nwnu.edu.cn (Z.P. Wang).}

\date{}\maketitle

\noindent{\footnotesize {\bf Abstract}  Let $\mathcal{C}$ be a triangulated category with a proper class $\xi$ of triangles. Asadollahi and Salarian introduced and studied $\xi$-Gorenstein projective and $\xi$-Gorenstein injective objects, and developed Gorenstein homological algebra in $\mathcal{C}$. In this paper, we further study Gorenstein homological properties for a triangulated category. First, we discuss the stability of $\xi$-Gorenstein projective objects, and show that the subcategory $\mathcal{GP}(\xi)$ of all $\xi$-Gorenstein projective objects has a strong stability. That is, an iteration of the procedure used to define the $\xi$-Gorenstein projective objects yields exactly the $\xi$-Gorenstein projective objects. Second,  we give some equivalent characterizations for $\xi$-Gorenstein projective dimension of object in $\mathcal{C}$.

\vspace{0.2cm}
\noindent{\footnotesize {2010 {\bf{Mathematics Subject
Classification}}:}
16E05, 16E10, 18G35

\noindent{\footnotesize {{\bf{Key words}}:} stability; $\xi$-Gorenstein projective object; triangulated category

\section{Introduction}
Triangulated categories were introduced independently in algebraic geometry by Verdier in his th\'{e}se \cite{Ve77}, and in algebraic topology by Puppe \cite{Pu62} in the early sixies, which have by now become indispensable in many areas of mathematics such as algebraic geometry, stable homotopy theory and representation theory \cite{BBD82,Ha88, Ma83}. The basic properties of triangulated categories can be found in Neeman's book \cite{Ne01}.

Let $\mathcal{C}$ be a triangulated category with triangles $\Delta$. Beligiannis \cite{Be00} developed homological algebra in $\mathcal{C}$ which parallels the homological algebra in an exact category in the sense of Quillen. By specifying a class of triangles $\xi\subseteq\Delta$, which is called a proper class of triangles, he introduced $\xi$-projective objects, $\xi$-projective dimensions and their duals.

Auslander and Bridger generalized in \cite{AB69} finitely generated projective modules to finitely generated modules of Gorenstein dimension zero over commutative noetherian rings. Furthermore, Enochs and Jenda introduced in \cite{EJ95} Gorenstein projective modules for arbitrary modules over a general ring, which is a generalization of finitely generated modules of Gorenstein dimension zero, and dually they defined Gorenstein injective modules. Gorenstein homological algebra has been extensively studied by many authors, see for example \cite{AM02,Ch00,EJ00,Ho04}.

As a natural generalization of modules of Gorenstein dimension zero, Beligiannis \cite{Bel00} defined the concept of an $\mathcal{X}$-Gorenstein object in an additive category $\mathcal{C}$ for a contravariantly finite subcategory $\mathcal{X}$ of $\mathcal{C}$ such that any $\mathcal{X}$-epic has kernal in $\mathcal{C}$. In an attempt to extend the theory, Asadollahi and Salarian \cite{AS04} introduced and studied $\xi$-Gorenstein projective and $\xi$-Gorenstein injective objects, and then $\xi$-Gorenstein projective and $\xi$-Gorenstein injective dimensions of objects in a triangulated category which are defined by modifying what Enochs and Jenda have done in an abelian category \cite{EJ00}.

Let $\mathcal{A}$ be an abelian category and $\mathcal{X}$ an additive full subcategory of $\mathcal{A}$. Sather-Wagstaff, Sharif and White introduced in \cite{SSW08} the Gorenstein category $\mathcal{G}(\mathcal{X})$, which is defined as $\mathcal{G}(\mathcal{X})=\{A ~\text{is an object in}~ \mathcal{A}~| ~\\ \text{there exists an exact sequence}~ \cdots\rightarrow X_{1}\rightarrow X_{0}\rightarrow X^{0}\rightarrow X^{1}\rightarrow \cdots \text{in }~ \mathcal{X}, ~\text{which is both}~ \Hom_{\mathcal{A}}(\mathcal{X},-)\text{-exact}\\ \text{and}~ \Hom_{\mathcal{A}}(-,\mathcal{X})~\text{-exact, such that}~ A\cong \im(X_{0}\rightarrow X^{0})\}$. Set $\mathcal{G}^{0}(\mathcal{X})=\mathcal{X}$, $\mathcal{G}^{1}(\mathcal{X})=\mathcal{G}(\mathcal{X})$, and inductively set $\mathcal{G}^{n+1}(\mathcal{X})=\mathcal{G}^{n}(\mathcal{G}(\mathcal{X}))$ for any $n\geq1$. They proved that when $\mathcal{X}$ is self-orthogonal, $\mathcal{G}^{n}(\mathcal{X})=\mathcal{G}(\mathcal{X})$ for any $n\geq1$; and they proposed the question whether $\mathcal{G}^{2}(\mathcal{X})=\mathcal{G}(\mathcal{X})$ holds for an arbitrary subcategory $\mathcal{X}$. See \cite[4.10, 5.8]{SSW08}. Recently, Huang \cite{Hu13} proved that the answer to this question is affirmative. This shows that $\mathcal{G}(\mathcal{X})$, in particular the subcategory $\mathcal{GP}(\mathcal{A})$ of all Gorenstein projective objects, has a strong stability. Kong and Zhang give a slight generalization of this stability by a different method \cite{KZ14}.

Inspired by the above results, we consider the stability of the subcategory $\mathcal{GP}(\xi)$ of all $\xi$-Gorenstein projective objects, which is introduced by Asadollahi and Salarian \cite{AS04}. Set $\mathcal{G}^{0}\mathcal{P}(\xi)=\mathcal{P}(\xi)$, $\mathcal{G}^{1}\mathcal{P}(\xi)=\mathcal{GP}(\xi)$, and inductively set $\mathcal{G}^{n+1}\mathcal{P}(\xi)=\mathcal{G}^{n}
(\mathcal{GP}(\xi))$ for any $n\geq1$. A natural question is whether $\mathcal{G}^{n}\mathcal{P}(\xi)
=\mathcal{GP}(\xi)$.

In Section $2$, we give some terminologies and some preliminary results. Section $3$ is devoted to answer the above question. We will prove the following theorem.

\noindent{\bf Main theorem} Let $\mathcal{C}$ be a triangulated category with enough $\xi$-projectives, where $\xi$ is a proper class of triangles. Then $\mathcal{G}^{n}\mathcal{P}(\xi)=\mathcal{GP}(\xi)$ for any $n\geq 1$.

The above theorem shows that the subcategory $\mathcal{GP}(\xi)$ of all $\xi$-Gorenstein projective objects has a strong stability. That is, an iteration of the procedure used to define the $\xi$-Gorenstein projective objects yields exactly the $\xi$-Gorenstein projective objects. Finally, we give some equivalent characterizations for $\xi$-Gorenstein projective dimension of an object $A$ in $\mathcal{C}$.

\section{Some basic facts in triangulated categories}
This section is devoted to recall the definitions and elementary properties of triangulated
categories used throughout the paper. For the terminology we shall follow \cite{ AS04, AS06,Be00}.

Let $\mathcal{C}$ be an additive category and $\Sigma: \mathcal{C}\rightarrow \mathcal{C}$ be an additive functor. Let $\Diag(\mathcal{C}, \Sigma)$ denote the category whose objects are diagrams in $\mathcal{C}$ of the form $A\rightarrow B\rightarrow C\rightarrow \Sigma A$, and morphisms between two objects $A_{i}\rightarrow B_{i}\rightarrow C_{i}\rightarrow \Sigma A_{i}$, $i=1,2$, are a triple of morphisms $\alpha:A_{1}\rightarrow A_{2}$, $\beta:B_{1}\rightarrow B_{2}$ and $\gamma:C_{1}\rightarrow C_{2}$, such that the following diagram commutes:
\begin{center}$\xymatrix{A_{1}\ar[r]^{f_{1}}\ar[d]_{\alpha}&
B_{1}\ar[r]^{g_{1}}\ar[d]_{\beta}&
C_{1}\ar[r]^{h_{1}}\ar[d]_{\gamma}&\Sigma A_{1}\ar[d]_{\Sigma\alpha}\\
A_{2}\ar[r]^{f_{2}}&B_{2}\ar[r]^{g_{2}}&C_{2}\ar[r]^{h_{2}}&\Sigma A_{2}.}$
\end{center}

A triangulated category is a triple $(\mathcal{C}, \Sigma, \Delta)$, where $\mathcal{C}$ is an additive category, $\Sigma$ is an autoequivalence of $\mathcal{C}$ and $\Delta$ is a full subcategory of $\Diag(\mathcal{C}, \Sigma)$ which satisfies the following axioms. The elements of $\Delta$ are then called triangles.

(Tr1) Every diagram isomorphic to a triangle is a triangle. Every morphism $f: A\rightarrow B$ in $\mathcal{C}$ can be embedded into a triangle $A\stackrel{f}\rightarrow B\stackrel{g}\rightarrow C\stackrel{h}\rightarrow \Sigma A$. For any object $A\in \mathcal{C}$, the diagram $A\stackrel{1_{A}}\rightarrow A\rightarrow 0\rightarrow \Sigma A$ is a triangle, where $1_{A}$ denotes the identity morphism from $A$ to $A$.

(Tr2) $A\stackrel{f}\rightarrow B\stackrel{g}\rightarrow C\stackrel{h}\rightarrow \Sigma A$ is a triangle if and only if $ B\stackrel{g}\rightarrow C\stackrel{h}\rightarrow \Sigma A\stackrel{-\Sigma f}\rightarrow \Sigma B$.

(Tr3) Given triangles $A_{i}\stackrel{f_{i}}\rightarrow B_{i}\stackrel{g_{i}}\rightarrow C_{i}\stackrel{h_{i}}\rightarrow \Sigma A_{i}$, $i=1,2$, and morphism $\alpha: A_{1}\rightarrow A_{2}$, $\beta: B_{1}\rightarrow B_{2}$ such that $f_{2}\alpha=\beta f_{1}$, there exists a morphism $\gamma: C_{1}\rightarrow C_{2}$ such that $(\alpha, \beta, \gamma)$ is a morphism from the first triangle to the second.

(Tr4) (The Octahedral Axiom) Given triangles $$A\stackrel{f}\rightarrow B\stackrel{i}\rightarrow C'\stackrel{i'}\rightarrow \Sigma A,~~~~B\stackrel{g}\rightarrow C\stackrel{j}\rightarrow A'\stackrel{j'}\rightarrow \Sigma B,~~~~A\stackrel{gf}\rightarrow C\stackrel{k}\rightarrow B'\stackrel{k'}\rightarrow \Sigma A,$$
there exist morphisms $f': C'\rightarrow B'$ and $g': B'\rightarrow A'$ such that the following diagram commutes and the third column is a triangle:
\begin{center}$\xymatrix{A\ar@{=}[d]\ar[r]^{f}&
B\ar[d]^{g}\ar[r]^{i}&C'\ar@{.>}[d]^{f'}\ar[r]^{i'}&\Sigma A\ar@{=}[d]\\
A\ar[r]^{gf}&
C\ar[d]^{j}\ar[r]^{k}&B'\ar@{.>}[d]^{g'}\ar[r]^{k'}&\Sigma A\ar[d]^{\Sigma f}\\
&
A'\ar[d]^{j'}\ar@{=}[r]&A'\ar[d]^{j'\Sigma i}\ar[r]^{j'}&\Sigma B\\
&\Sigma B\ar[r]^{\Sigma i}&\Sigma C'&}$
\end{center}

 Some equivalent formulations for the Octahedral Axiom $(\mathrm{Tr4})$ are given in \cite[2.1]{Be00}, which are more convient to use. If $\mathcal{C}=(\mathcal{C}, \Sigma, \Delta)$ satisfies all the axioms of a triangulated category except possible of $(\mathrm{Tr4})$, then $(\mathrm{Tr4})$ is equivalent to each of the following:

 $\mathbf{Base~~Change}$: For any triangle $A\stackrel{f}\rightarrow
 B\stackrel{g}\rightarrow C\stackrel{h}\rightarrow  \Sigma A \in \Delta$ and morphism $\varepsilon: E\rightarrow C$, there exists a commutative diagram:
 \begin{center}
$\xymatrix{
    & M \ar[d]_{\alpha}\ar@{=}[r] &M \ar[d]_{\delta} &\\
      A\ar@{=}[d]_{} \ar[r]^{f'} & G \ar[d]_{\beta}\ar[r]^{g'} &E \ar[d]_{\varepsilon} \ar[r]^{h'}&\Sigma A \ar@{=}[d]\\
     A \ar[r]_{f} & B\ar[d]_{\gamma}\ar[r]_{g} &C \ar[d]_{\zeta} \ar[r]_{h}&\Sigma A \ar[d]_{\Sigma f'} \\
       & \Sigma M \ar@{=}[r] &\Sigma M \ar[r]_{-\Sigma\alpha}&\Sigma G}$
\end{center}
in which all horizontal and vertical diagrams are triangles in $\Delta$.

  $\mathbf{Cobase~~Change}$: For any triangle $A\stackrel{f}\rightarrow
 B\stackrel{g}\rightarrow C\stackrel{h}\rightarrow  \Sigma A \in \Delta$ and morphism $\alpha: A\rightarrow D$, there exists a commutative diagram:
 \begin{center}
$\xymatrix{
    & N \ar[d]_{\zeta}\ar@{=}[r] &N \ar[d]_{\delta} &\\
      \Sigma^{-1}C\ar@{=}[d]_{} \ar[r]^{-\Sigma^{-1}h} & A \ar[d]_{\alpha}\ar[r]^{f} &B\ar[d]_{\beta} \ar[r]^{g}&C\ar@{=}[d]\\
     \Sigma^{-1}C \ar[r]_{-\Sigma^{-1}h'} & D\ar[d]_{\eta}\ar[r]_{f'} &F \ar[d]_{\upsilon} \ar[r]_{g'}&C \ar[d]_{-h} \\
       & \Sigma N \ar@{=}[r] &\Sigma N \ar[r]_{-
       \Sigma \eta}&\Sigma A}$
\end{center}
in which all horizontal and vertical diagrams are triangles in $\Delta$.

The following definitions are quoted verbatim from \cite[Section 2]{Be00}. A class of triangles $\xi$ is closed under base change if for any triangle $A\stackrel{f}\rightarrow
 B\stackrel{g}\rightarrow C\stackrel{h}\rightarrow  \Sigma A \in \xi$ and any morphism $\varepsilon: E\rightarrow C$ as in the base change diagram above, the triangle $A\stackrel{f'}\rightarrow
 G\stackrel{g'}\rightarrow E\stackrel{h'}\rightarrow  \Sigma A $ belongs to $\xi$. Dually, a class of triangles $\xi$ is closed under cobase change if for any triangle $A\stackrel{f}\rightarrow
 B\stackrel{g}\rightarrow C\stackrel{h}\rightarrow  \Sigma A \in \xi$ and any morphism $\alpha: A\rightarrow D$ as in the cobase change diagram above, the triangle $D\stackrel{f'}\rightarrow
 F\stackrel{g'}\rightarrow C\stackrel{h'}\rightarrow  \Sigma D $ belongs to $\xi$. A class of triangles $\xi$ is closed under suspension if for any triangle $A\stackrel{f}\rightarrow
 B\stackrel{g}\rightarrow C\stackrel{h}\rightarrow  \Sigma A \in \xi$ and $i\in \mathbb{Z}$, the triangle $\Sigma^{i}A\xrightarrow{(-1)^{i}\Sigma^{i}f}
 \Sigma^{i}B\xrightarrow {(-1)^{i}\Sigma^{i}g} \Sigma^{i}C\xrightarrow{(-1)^{i}\Sigma^{i}h}  \Sigma^{i+1} A$ is in $\xi$. A class of triangles $\xi$ is called saturated if in the situation of base change, whenever the third vertical and the second horizontal triangles are in $\xi$, then the triangle $A\stackrel{f}\rightarrow
 B\stackrel{g}\rightarrow C\stackrel{h}\rightarrow  \Sigma A $ is in $\xi$. An easy consequence of the octahedral axiom is that $\xi$ is saturated if and only if in the situation of cobase change, whenever the second vertical and the third horizontal triangles are in $\xi$, then the triangle $A\stackrel{f}\rightarrow
 B\stackrel{g}\rightarrow C\stackrel{h}\rightarrow  \Sigma A $ is in $\xi$.

 A triangle $(T): A\stackrel{f}\rightarrow
 B\stackrel{g}\rightarrow C\stackrel{h}\rightarrow  \Sigma A \in \Delta$ is called split if it is isomorphic to the triangle $A\stackrel{(1~0)}\longrightarrow
 A\oplus C\stackrel{{0\choose 1}}\longrightarrow C\stackrel{0}\rightarrow  \Sigma A$. It is easy to see that $(T)$ is split if and only if $f$ is a section or $g$ is a retraction or $h=0$. The full subcategory of $\Delta$ consisting of the split triangles will be denoted by $\Delta_{0}$.

 \begin{definition}(\cite[2.2]{Be00}) Let $\mathcal{C}=(\mathcal{C}, \Sigma, \Delta)$ be a triangulated category. A class $\xi\subseteq\Delta$ is called a proper class of triangles if the following conditions hold.

 (1) $\xi$ is closed under isomorphisms, finite coproducts and $\Delta_{0}\subseteq\xi\subseteq\Delta$.

 (2) $\xi$ is closed under suspensions and is saturated.

 (3) $\xi$ is closed under base and cobase change.

 \end{definition}

 It is known that $\Delta_{0}$ and the class of all triangles $\Delta$ in $\mathcal{C}$ are proper classes of triangles. There are more interesting examples of proper classes of triangles enumerated in \cite[Example 2.3]{Be00}.

 Throughout we fix a proper class $\xi$ of triangles in the triangulated category $\mathcal{C}$.

 \begin{definition}(\cite[4.1]{Be00}) An object $P\in \mathcal{C}$ (resp., $I\in \mathcal{C}$) is called $\xi$-projective (resp., $\xi$-injective) if for any triangle $A\rightarrow
 B\rightarrow C\rightarrow  \Sigma A$ in $\xi$, the induced sequence of abelian groups $0\rightarrow \Hom_{\mathcal{C}}(P,A)\rightarrow\Hom_{\mathcal{C}}(P,B)
 \rightarrow\Hom_{\mathcal{C}}(P,C)\rightarrow0$ (resp., $0\rightarrow \Hom_{\mathcal{C}}(C,I)\rightarrow\Hom_{\mathcal{C}}(B,I)
 \rightarrow\Hom_{\mathcal{C}}(A,I)\rightarrow0$) is exact.
 \end{definition}

 It follows easily from the definition that the subcategory $\mathcal{P}(\xi)$ of $\xi$-projective objects and the subcategory $\mathcal{I}(\xi)$ of $\xi$-injective objects are full, additive, $\Sigma$-stable, and closed under isomorphisms and direct summands. The category $\mathcal{C}$ is said to have enough $\xi$-projectives (resp., $\xi$-injectives) if for any object $A\in \mathcal{C}$, there exists a triangle $K\rightarrow
 P\rightarrow A\rightarrow  \Sigma K$ (resp., $A\rightarrow
 I\rightarrow L\rightarrow  \Sigma A$) in $\xi$ with $P\in \mathcal{P}(\xi)$ (resp., $I\in \mathcal{I}(\xi)$). However, it is not so easy to find a proper class $\xi$ of triangles in a triangulated category having enough $\xi$-projectives or $\xi$-injectives in general. Here we present a nontrivial example due to Beligiannis, and also, see more nontrivial examples in \cite[Sections 12.4 and 12.5]{Be00}, which are of great interest. Take $\mathcal{C}$ to be the unbounded homotopy category of complexes of objects from a Grothendieck category which has enough projectives. Then the so-called Cartan-Eilenberg projective and injective complexes form the relative projective and injective objects for a proper class of triangles in $\mathcal{C}$.

 The following lemma is quoted from \cite[4.2]{Be00}.

 \begin{lemma}\label{lem1} Let $\mathcal{C}$ have enough $\xi$-projectives. Then a triangle $A\rightarrow
 B\rightarrow C\rightarrow  \Sigma A$ is in $\xi$ if and only if for all $P\in \mathcal{P}(\xi)$ the induced sequence $0\rightarrow \Hom_{\mathcal{C}}(P,A)\rightarrow\Hom_{\mathcal{C}}(P,B)
 \rightarrow\Hom_{\mathcal{C}}(P,C)\rightarrow0$ is exact.
 \end{lemma}

 Recall that an $\xi$-exact sequence $X$ is a diagram
 $$\cdots\longrightarrow X_{1}\stackrel{d_{1}}\longrightarrow X_{0}\stackrel{d_{0}}\longrightarrow X_{-1}\stackrel{d_{-1}}\longrightarrow X_{-2}\longrightarrow\cdots$$
 in $\mathcal{C}$, such that for each $n\in \mathbb{Z}$, there exists triangle $K_{n+1}\stackrel{f_{n}}\longrightarrow X_{n}\stackrel{g_{n}}\longrightarrow K_{n}\stackrel{h_{n}}\longrightarrow \Sigma K_{n+1}$ in $\xi$ and the differential is defined as $d_{n}=f_{n-1}g_{n}$ for any $n$.

 Let $\mathcal{X}$ be a full and additive subcategory of $\mathcal{C}$. A triangle $A\rightarrow
 B\rightarrow C\rightarrow  \Sigma A$ in $\xi$ is called $\Hom_{\mathcal{C}}(-,\mathcal{X})$-exact (resp., $\Hom_{\mathcal{C}}(\mathcal{X},-)$-exact), if for any $X\in \mathcal{X}$, the induced sequence of abelian groups  $0\rightarrow \Hom_{\mathcal{C}}(C, X)\rightarrow\Hom_{\mathcal{C}}(B, X)
 \rightarrow\Hom_{\mathcal{C}}(A, X)\rightarrow0$ (resp., $0\rightarrow \Hom_{\mathcal{C}}(X, A)\rightarrow\Hom_{\mathcal{C}}(X,B)
 \rightarrow\Hom_{\mathcal{C}}(X,C)\rightarrow0$) is exact.

 A complete $\xi$-projective resolution is a diagram
 $$\mathbf{P}:\cdots \rightarrow P_{1}\stackrel{d_{1}}\longrightarrow P_{0}
 \stackrel{d_{0}}\longrightarrow P_{-1}\rightarrow \cdots$$
 in $\mathcal{C}$ such that for any integer $n$, $P_{n}\in \mathcal{P}(\xi)$ and there exist $\Hom_{\mathcal{C}}(-,\mathcal{P}(\xi))$-exact triangles
 $$K_{n+1}\stackrel{f_{n}}\longrightarrow P_{n}\stackrel{g_{n}}\longrightarrow K_{n}\stackrel{h_{n}}\longrightarrow \Sigma K_{n+1}$$ in $\xi$ and the differential is defined as $d_{n}=f_{n-1}g_{n}$ for any $n$.
 \begin{definition}(see \cite{AS04}) Let $\mathbf{P}$ be a complete $\xi$-projective resolution in $\mathcal{C}$. So for any integer $n$, there exist triangles
 $$K_{n+1}\stackrel{f_{n}}\longrightarrow P_{n}\stackrel{g_{n}}\longrightarrow K_{n}\stackrel{h_{n}}\longrightarrow \Sigma K_{n+1}$$
in $\xi$. The objects $K_{n}$ for any integer $n$, are called $\xi$-Gorenstein projective ($\xi$-$\mathcal{G}$projective for short).
 \end{definition}
Dually, one can define complete $\xi$-injective resolution and $\xi$-Gorenstein injective ($\xi$-$\mathcal{G}$injective for short) objects.

We denote by $\mathcal{GP}(\xi)$ and $\mathcal{GI}(\xi)$ the subcategory of $\xi$-$\mathcal{G}$projective and $\xi$-$\mathcal{G}$injective objects of $\mathcal{C}$ respectively. It is obvious that $\mathcal{P}(\xi)\subseteq\mathcal{GP}(\xi)$ and $\mathcal{I}(\xi)\subseteq\mathcal{GI}(\xi)$; $\mathcal{GP}(\xi)$ and $\mathcal{GI}(\xi)$ are full, additive, $\Sigma$-stable, and closed under isomorphisms and direct summands.

\section{Main results}

Throughout the paper, $\mathcal{C}$ is a triangulated category with enough $\xi$-projectives and enough $\xi$-injectives, where $\xi$ is a fixed proper class of triangles.

Recall that in \cite{AS04} a diagram $$\mathbf{P}:\cdots \rightarrow P_{1}\stackrel{d_{1}}\longrightarrow P_{0}
 \stackrel{d_{0}}\longrightarrow P_{-1}\rightarrow \cdots$$
 in $\mathcal{C}$ is a complete $\xi$-projective resolution, if for any $n$, $P_{n}\in \mathcal{P}(\xi)$, and there exist $\Hom_{\mathcal{C}}(-,\mathcal{P}(\xi))$-exact triangles
 $$K_{n+1}\stackrel{f_{n}}\longrightarrow P_{n}\stackrel{g_{n}}\longrightarrow K_{n}\stackrel{h_{n}}\longrightarrow \Sigma K_{n+1}$$ in $\xi$ and the differential is defined as $d_{n}=f_{n-1}g_{n}$ for any $n$. That is, $\mathbf{P}$ is an $\xi$-exact sequence of $\xi$-projective objects, and is $\Hom_{\mathcal{C}}(-,\mathcal{P}(\xi))$-exact.
  The objects $K_{n}$ for any integer $n$, are called $\xi$-Gorenstein projective ($\xi$-$\mathcal{G}$projective for short). Dually, one can define complete $\xi$-injective resolution and $\xi$-Gorenstein injective ($\xi$-$\mathcal{G}$injective for short) objects.

We study only the case of $\xi$-$\mathcal{G}$projective objects since the study of the $\xi$-$\mathcal{G}$injective objects is dual.

An $\xi$-$\mathcal{G}$projective resolution of $A\in \mathcal{C}$ is an $\xi$-exact sequence $\cdots \rightarrow G_{n}\rightarrow G_{n-1}\rightarrow \cdots \rightarrow G_{1}\rightarrow G_{0}\rightarrow A\rightarrow 0$ in $\mathcal{C}$, such that $G_{n}\in \mathcal{GP}(\xi)$ for all $n\geq 0$. The definition is different from \cite[Definition 4.2]{AS04}.

\begin{lemma}\label{lem4} Let $\mathcal{C}$ be a triangulated category with enough $\xi$-projectives, $A\in \mathcal{C}$. Then $A$ has an $\xi$-projective resolution which is $\Hom_{\mathcal{C}}(-, \mathcal{P}(\xi))$-exact if and only if $A$ has an $\xi$-$\mathcal{G}$projective resolution which is $\Hom_{\mathcal{C}}(-, \mathcal{P}(\xi))$-exact.
\end{lemma}

\noindent{\bf Proof} Since $\mathcal{P}(\xi)\subseteq \mathcal{GP}(\xi)$, it is enough to show the ``if" part. Assume that $A$ has an $\xi$-$\mathcal{G}$projective resolution which is $\Hom_{\mathcal{C}}(-, \mathcal{P}(\xi))$-exact. Then there exists a triangle $\delta: B\rightarrow G_{0}\rightarrow A\rightarrow \Sigma B\in \xi$ which is $\Hom_{\mathcal{C}}(-, \mathcal{P}(\xi))$-exact, where $G_{0}\in \mathcal{GP}(\xi)$ and $B$ has an $\xi$-$\mathcal{G}$projective resolution which is $\Hom_{\mathcal{C}}(-, \mathcal{P}(\xi))$-exact. Since $G_{0}\in \mathcal{GP}(\xi)$, there exists a triangle $\eta: G_{0}'\rightarrow P_{0}\rightarrow G_{0}\rightarrow \Sigma G_{0}'\in \xi$ such that $G_{0}'\in \mathcal{GP}(\xi)$, $P_{0}\in \mathcal{P}(\xi)$ and it is $\Hom_{\mathcal{C}}(-, \mathcal{P}(\xi))$-exact. By base change, we have the following commutative diagram:
\begin{center}
$\xymatrix{&&  \Sigma^{-1}A\ar[d]_{}\ar@{=}[r] &\Sigma^{-1}A\ar[d]_{} &\\
   \eta':&   G_{0}'\ar@{=}[d]_{} \ar[r]^{} & L \ar[d]_{}\ar[r]^{} &B\ar[d]_{} \ar[r]^{}&\Sigma G_{0}'\ar@{=}[d]\\
    \eta:& G_{0}'\ar[r]_{} & P_{0}\ar[d]_{}\ar[r]^{} &G_{0} \ar[d]_{} \ar[r]_{}&\Sigma G_{0}'\\
       && \Sigma A\ar@{=}[r] &A &}$
\end{center}
Since $\xi$ is closed under base change, we have $\eta'\in \xi$. Applying the functor $\Hom_{\mathcal{C}}(\mathcal{P}(\xi), -)$ to the above diagram, we obtain the following commutative diagram
\begin{center}
$\xymatrix{
&&&0\ar@{.>}[d]&0\ar[d]\\
  \Hom_{\mathcal{C}}(\mathcal{P}(\xi), \eta'):& 0\ar@{.>}[r] & \Hom_{\mathcal{C}}(\mathcal{P}(\xi), G_{0}')\ar@{=}[d]_{}\ar[r] &\Hom_{\mathcal{C}}(\mathcal{P}(\xi), L) \ar[d]_{} \ar[r]&\Hom_{\mathcal{C}}(\mathcal{P}(\xi), B)\ar[d]\ar@{.>}[r]&0\\
 \Hom_{\mathcal{C}}(\mathcal{P}(\xi), \eta):& 0\ar[r] & \Hom_{\mathcal{C}}(\mathcal{P}(\xi), G_{0}')\ar[r] &\Hom_{\mathcal{C}}(\mathcal{P}(\xi), P_{0}) \ar[d]_{} \ar[r]&\Hom_{\mathcal{C}}(\mathcal{P}(\xi), G_{0})\ar[d]\ar[r]&0\\
       &&& \Hom_{\mathcal{C}}(\mathcal{P}(\xi), A) \ar@{.>}[d]\ar@{=}[r] &\Hom_{\mathcal{C}}(\mathcal{P}(\xi), A) \ar[d]&\\
       &&&0&0}$
\end{center}
By snake lemma, we get that the second vertical sequence is exact, and so $L\rightarrow P_{0}\rightarrow A\rightarrow \Sigma L\in \xi$ by Lemma \ref{lem1}. For any $Q\in \mathcal{P}(\xi)$, applying $\Hom_{\mathcal{C}}(-, Q)$ to the above base change diagram, we obtain the following commutative diagram:
\begin{center}
$\xymatrix{&&&\Hom_{\mathcal{C}}(\eta, Q)&\Hom_{\mathcal{C}}(\eta', Q)\\
&&&0\ar[d]&0\ar@{.>}[d]\\
 \Hom_{\mathcal{C}}(\delta, Q):&  0\ar[r] & \Hom_{\mathcal{C}}(A, Q)\ar@{=}[d]_{}\ar[r] &\Hom_{\mathcal{C}}(G_{0}, Q) \ar[d]_{} \ar[r]&\Hom_{\mathcal{C}}(B, Q)\ar[d]\ar[r]&0\\
  \Hom_{\mathcal{C}}(\delta', Q):& 0\ar@{.>}[r] & \Hom_{\mathcal{C}}(A, Q)\ar[r] &\Hom_{\mathcal{C}}(P_{0}, Q) \ar[d]_{} \ar[r]&\Hom_{\mathcal{C}}(L, Q)\ar[d]\ar@{.>}[r]&0\\
      & && \Hom_{\mathcal{C}}(G_{0}', Q) \ar[d]\ar@{=}[r] &\Hom_{\mathcal{C}}(G_{0}', Q) \ar@{.>}[d]&\\
       &&&0&0}$
\end{center}
By snake lemma, one can get that $\Hom_{\mathcal{C}}(\delta', Q)$ and $\Hom_{\mathcal{C}}(\eta', Q)$ are exact. Since $B$ has an $\xi$-$\mathcal{G}$projective resolution which is $\Hom_{\mathcal{C}}(-, \mathcal{P}(\xi))$-exact, there exists a triangle $C\rightarrow G_{1}\rightarrow B\rightarrow \Sigma C\in \xi$ which is $\Hom_{\mathcal{C}}(-, \mathcal{P}(\xi))$-exact, where $G_{1}\in \mathcal{GP}(\xi)$ and $C$ has an $\xi$-$\mathcal{G}$projective resolution which is $\Hom_{\mathcal{C}}(-, \mathcal{P}(\xi))$-exact. By base change and the similar method above, we have the following commutative diagram:
\begin{center}
$\xymatrix{&&\zeta'&\zeta&\\
    &&  C\ar[d]_{}\ar@{=}[r] &C \ar[d]_{} &\\
   \epsilon':&   G_{0}'\ar@{=}[d]_{} \ar[r]^{} & M \ar[d]_{}\ar[r]^{} &G_{1}\ar[d]_{} \ar[r]^{}& \Sigma G_{0}'\ar@{=}[d]\\
    \epsilon':& G_{0}' \ar[r]_{} & L\ar[d]_{}\ar[r]^{} &B \ar[d]_{} \ar[r]_{}& \Sigma G_{0}'\\
       && \Sigma C \ar@{=}[r] &\Sigma C &}$
\end{center}
such that $\epsilon'$ and $\zeta'$ are in $\xi$. Since $G_{0}'\in \mathcal{GP}(\xi) $ and $G_{1}\in \mathcal{GP}(\xi)$, by \cite[Theorem 3.11]{AS04}, we have $M\in \mathcal{GP}(\xi)$. Thus $L$ has an $\xi$-$\mathcal{G}$projective resolution which is $\Hom_{\mathcal{C}}(-, \mathcal{P}(\xi))$-exact. Note that the triangle $L\rightarrow P_{0}\rightarrow A\rightarrow \Sigma L$ is $\Hom_{\mathcal{C}}(-, \mathcal{P}(\xi))$-exact. By repeating the preceding process, we have that $A$ has an $\xi$-projective resolution which is $\Hom_{\mathcal{C}}(-, \mathcal{P}(\xi))$-exact, as required. \hfill$\Box$

An $\xi$-projective ($\xi$-$\mathcal{G}$projective) coresolution of $A\in \mathcal{C}$ is an $\xi$-exact sequence
$0\rightarrow A\rightarrow X^{0}\rightarrow X^{1}\rightarrow \cdots $ in $\mathcal{C}$, such that $X^{n} \in \mathcal{P}(\xi)$ ($X^{n} \in\mathcal{GP}(\xi)$) for all $n\geq 0$.

\begin{lemma}\label{lem5} Let $\mathcal{C}$ be a triangulated category with enough $\xi$-projectives, $A\in \mathcal{C}$.  Then $A$ has an $\xi$-projective coresolution which is $\Hom_{\mathcal{C}}(-, \mathcal{P}(\xi))$-exact if and only if $A$ has an $\xi$-$\mathcal{G}$projective coresolution which is $\Hom_{\mathcal{C}}(-, \mathcal{P}(\xi))$-exact.
\end{lemma}

\noindent{\bf Proof} It is completely dual to the proof of Lemma \ref{lem4}. So we omit it.\hfill$\Box$

\begin{lemma} Let $\mathcal{C}$ be a triangulated category with enough $\xi$-projectives, $A\in \mathcal{C}$. Then the following statements are equivalent:

 (1) $A$ is an $\xi$-$\mathcal{G}$projective object.

  (2) $A$ has an $\xi$-projective resolution which is $\Hom_{\mathcal{C}}(-, \mathcal{P}(\xi))$-exact and has an $\xi$-projective coresolution which is $\Hom_{\mathcal{C}}(-, \mathcal{P}(\xi))$-exact.

  (3) There exist  $\Hom_{\mathcal{C}}(-, \mathcal{P}(\xi))$-exact triangles $K_{n+1}\rightarrow P_{n}\rightarrow K_{n}\rightarrow \Sigma K_{n+1} \in \xi$ such that $P_{n}\in \mathcal{P}(\xi)$ and $K_{0}=A$.
\end{lemma}

\noindent{\bf Proof} It follows from the definition of $\xi$-$\mathcal{G}$projective object. \hfill$\Box$

\begin{theorem} Let $\mathcal{C}$ be a triangulated category with enough $\xi$-projectives, $A\in \mathcal{C}$. The following are equivalent:

(1) $A$ is $\xi$-$\mathcal{G}$projective.

(2)  There exist  $\Hom_{\mathcal{C}}(-, \mathcal{GP}(\xi))$-exact and $\Hom_{\mathcal{C}}(\mathcal{GP}(\xi),-)$-exact triangles $K_{n+1}\rightarrow G_{n}\rightarrow K_{n}\rightarrow \Sigma K_{n+1} \in \xi$ such that $G_{n}\in \mathcal{GP}(\xi)$ and $K_{0}=A$.

(3) There exist  $\Hom_{\mathcal{C}}(-, \mathcal{GP}(\xi))$-exact triangles $K_{n+1}\rightarrow G_{n}\rightarrow K_{n}\rightarrow \Sigma K_{n+1} \in \xi$ such that $G_{n}\in \mathcal{GP}(\xi)$ and $K_{0}=A$.

(4) There exist  $\Hom_{\mathcal{C}}(-, \mathcal{P}(\xi))$-exact triangles $K_{n+1}\rightarrow G_{n}\rightarrow K_{n}\rightarrow \Sigma K_{n+1} \in \xi$ such that $G_{n}\in \mathcal{GP}(\xi)$ and $K_{0}=A$.

(5) There exists a $\Hom_{\mathcal{C}}(\mathcal{GP}(\xi),-)$-exact and $\Hom_{\mathcal{C}}(-, \mathcal{GP}(\xi))$-exact triangle $A\rightarrow G\rightarrow A\rightarrow \Sigma A \in \xi$ such that $G\in \mathcal{GP}(\xi)$.

(6) There exists a $\Hom_{\mathcal{C}}(-, \mathcal{GP}(\xi))$-exact triangle $A\rightarrow G\rightarrow A\rightarrow \Sigma A \in \xi$ such that $G\in \mathcal{GP}(\xi)$.

(7) There exists a $\Hom_{\mathcal{C}}(-, \mathcal{P}(\xi))$-exact triangle $A\rightarrow G\rightarrow A\rightarrow \Sigma A \in \xi$ such that $G\in \mathcal{GP}(\xi)$.
\end{theorem}

\noindent{\bf Proof}

(1)$\Rightarrow$(2) Let $A$ be an $\xi$-$\mathcal{G}$projective object of $\mathcal{C}$. Consider the triangles
\begin{center}$0\rightarrow A\stackrel{1}\rightarrow A\rightarrow 0$ and $A\stackrel{1}\rightarrow A\rightarrow 0\rightarrow \Sigma A$
\end{center}
Since $\xi$ is proper, it contains $\Delta_{0}$. So the above two triangles are in $\xi$. It is easy to see that they are $\Hom_{\mathcal{C}}(-, \mathcal{GP}(\xi))$-exact and $\Hom_{\mathcal{C}}(\mathcal{GP}(\xi),-)$-exact triangles.

(2)$\Rightarrow$(3) and (3)$\Rightarrow$(4) are clear.

(4)$\Rightarrow$(1) It follows from Lemma \ref{lem4} and Lemma \ref{lem5}.

(1)$\Rightarrow$(5) Let $A$ be an $\xi$-$\mathcal{G}$projective object of $\mathcal{C}$. Consider the split triangle
\begin{center}$\delta:~~A\stackrel{{1\choose 0}}\rightarrow A\oplus A\stackrel{(0 ~1)}\rightarrow A\stackrel{0}\rightarrow \Sigma A$
\end{center}
Since $\xi$ is proper, it contains $\Delta_{0}$. So $\delta$ is in $\xi$. For any $Q\in \mathcal{GP}(\xi)$, applying the functors  $\Hom_{\mathcal{C}}(Q, -)$ and $\Hom_{\mathcal{C}}(-, Q)$ to the above triangle, we get the following exact sequences $$0\rightarrow\Hom_{\mathcal{C}}(A,Q)
\rightarrow\Hom_{\mathcal{C}}(A\oplus A,Q)\rightarrow\Hom_{\mathcal{C}}(A,Q)\rightarrow0,$$
 $$0\rightarrow\Hom_{\mathcal{C}}(Q,A)
\rightarrow\Hom_{\mathcal{C}}(Q, A\oplus A)\rightarrow\Hom_{\mathcal{C}}(Q,A)\rightarrow0.$$ By \cite[Theorem 3.11]{AS04}, we can obtain $A\oplus A \in \mathcal{GP}(\xi)$. So we are done.

(5)$\Rightarrow$(6)$\Rightarrow$(7)$\Rightarrow$(1) are clear.\hfill$\Box$

Denote by $\mathcal{GP}(\xi)$ the subcategory of all $\xi$-$\mathcal{G}$projective objects. Set $\mathcal{G}^{0}\mathcal{P}(\xi)=\mathcal{P}(\xi)$, $\mathcal{G}^{1}\mathcal{P}(\xi)=\mathcal{GP}(\xi)$, and inductively set $\mathcal{G}^{n+1}\mathcal{P}(\xi)=\mathcal{G}^{n}
(\mathcal{GP}(\xi))$ for any $n\geq1$. Now we can obtain our main theorem.

\begin{theorem} \label{mainth} Let $\mathcal{C}$ be a triangulated category with enough $\xi$-projectives. Then $\mathcal{G}^{n}\mathcal{P}(\xi)=\mathcal{GP}(\xi)$ for any $n\geq 1$.
\end{theorem}

\noindent{\bf Proof} It is easy to see that $\mathcal{P}(\xi)\subseteq \mathcal{GP}(\xi) \subseteq \mathcal{G}^{2}\mathcal{P}(\xi)\subseteq \mathcal{G}^{3}\mathcal{P}(\xi)\subseteq \cdots$ is an ascending chain of subcategories of $\mathcal{C}$. By $(1)\Leftrightarrow (3)$ of the above theorem, we have that $\mathcal{G}^{2}\mathcal{P}(\xi)=\mathcal{GP}(\xi)$. By using induction on $n$ we get easily the assertion. \hfill$\Box$

In order to give some equivalent characterizations for $\xi$-Gorenstein projective dimension of an object $A$ in $\mathcal{C}$, one needs the following lemma.

\begin{lemma}\label{prop1} Let $0\rightarrow B\rightarrow G_{1}\rightarrow G_{0}\rightarrow A\rightarrow 0$ be an $\xi$-exact sequence with $G_{1},G_{0} \in \mathcal{GP}(\xi)$. Then there exist the following $\xi$-exact sequences:
$$0\rightarrow B\rightarrow P\rightarrow G_{0}'\rightarrow A\rightarrow 0$$
and $$0\rightarrow B\rightarrow G_{1}'\rightarrow Q\rightarrow A\rightarrow 0$$
where $P,Q\in \mathcal{P}(\xi)$ and $G_{0}',G_{1}'\in \mathcal{GP}(\xi).$
\end{lemma}

\noindent{\bf Proof} Since $G_{1}$ is in $\mathcal{GP}(\xi)$, there exists a triangle $G_{1}\rightarrow P\rightarrow G_{2}\rightarrow \Sigma G_{1} \in \xi$ with $P$ $\xi$-projective and $G_{2}$ $\xi$-$\mathcal{G}$projective. Since $0\rightarrow B\rightarrow G_{1}\rightarrow G_{0}\rightarrow A\rightarrow 0$ is an $\xi$-exact sequence, there exist triangles $B\rightarrow G_{1}\rightarrow K\rightarrow \Sigma B \in \xi$ and $K\rightarrow G_{0}\rightarrow A \rightarrow \Sigma K \in \xi$.
Then we have the following commutative diagrams by cobase change.
\begin{center}$\xymatrix{&B\ar@{=}[r]\ar[d]&B\ar[d]\\
\Sigma^{-1}G_{2}\ar[r]\ar@{=}[d]&G_{1}\ar[r]\ar[d]&P\ar[r]\ar[d]&
G_{2}\ar@{=}[d]\\
\Sigma^{-1}G_{2}\ar[r]&K\ar[r]\ar[d]&C\ar[r]\ar[d]&G_{2}\\
&\Sigma B\ar@{=}[r]&\Sigma B\\
}$
\end{center}

\begin{center}$\xymatrix{&\Sigma^{-1}G_{2}\ar@{=}[r]\ar[d]&
\Sigma^{-1}G_{2}\ar[d]\\
\Sigma^{-1}A\ar[r]\ar@{=}[d]&K\ar[r]\ar[d]&G_{0}\ar[r]\ar[d]&
A\ar@{=}[d]\\
\Sigma^{-1}A\ar[r]&C\ar[r]\ar[d]&G_{0}'\ar[r]\ar[d]&A\\
&G_{2}\ar@{=}[r]&G_{2}\\
}$
\end{center}
Since $\xi$ is closed under cobase change, we get $K\rightarrow C\rightarrow G_{2}\rightarrow \Sigma K \in \xi$ and $C\rightarrow G_{0}'\rightarrow A\rightarrow \Sigma C \in \xi$.

For any $Q\in \mathcal{P}(\xi)$, applying the functor $\Hom_{\mathcal{C}}(Q,-)$ to the above diagrams, we have the following commutative diagrams:
\begin{center}
$\xymatrix{
&0\ar[d]&0\ar@{.>}[d]&\\
& \Hom_{\mathcal{C}}(Q, B) \ar[d]\ar@{=}[r] &\Hom_{\mathcal{C}}(Q, B) \ar[d]&\\
0\ar[r] & \Hom_{\mathcal{C}}(Q, G_{1})\ar[d]_{}\ar[r] &\Hom_{\mathcal{C}}(Q, P) \ar[d]_{} \ar[r]&\Hom_{\mathcal{C}}(Q, G_{2})\ar@{=}[d]\ar[r]&0\\
0\ar@{.>}[r] & \Hom_{\mathcal{C}}(Q, K)\ar[r]\ar[d] &\Hom_{\mathcal{C}}(Q, C) \ar@{.>}[d] \ar[r]&\Hom_{\mathcal{C}}(Q, G_{2})\ar@{.>}[r]&0\\
&0&0&}$
\end{center}
\begin{center}
$\xymatrix{
&0\ar[d]&0\ar@{.>}[d]&\\
0\ar[r] & \Hom_{\mathcal{C}}(Q, K)\ar[d]_{}\ar[r] &\Hom_{\mathcal{C}}(Q, G_{0}) \ar[d]_{} \ar[r]&\Hom_{\mathcal{C}}(Q, A)\ar@{=}[d]\ar[r]&0\\
0\ar@{.>}[r] & \Hom_{\mathcal{C}}(Q, C)\ar[r]\ar[d] &\Hom_{\mathcal{C}}(Q, G_{0}') \ar[d] \ar[r]&\Hom_{\mathcal{C}}(Q, A)\ar@{.>}[r]&0\\& \Hom_{\mathcal{C}}(Q, G_{2}) \ar[d]\ar@{=}[r] &\Hom_{\mathcal{C}}(Q, G_{2}) \ar@{.>}[d]&\\
&0&0&}$
\end{center}
By snake lemma and Lemma \ref{lem1}, we have $B\rightarrow P\rightarrow C\rightarrow \Sigma B \in \xi$ and $G_{0}\rightarrow G_{0}'\rightarrow G_{2}\rightarrow \Sigma G_{0} \in \xi$.
Because both $G_{0}$ and $G_{2}$ are in $\mathcal{GP}(\xi)$, $G_{0}'$ is also in $\mathcal{GP}(\xi)$ by \cite[Theorem 3.11]{AS04}. Connecting the triangles $B\rightarrow P\rightarrow C\rightarrow \Sigma B$ and $C\rightarrow G_{0}'\rightarrow A\rightarrow \Sigma C$, we get the first desired $\xi$-exact sequence.

Since $G_{0}$ is $\xi$-$\mathcal{G}$projective, there is a triangle $G_{3}\rightarrow Q\rightarrow G_{0}\rightarrow \Sigma G_{3}\in \xi$ with $Q\in \mathcal{P}(\xi)$ and $G_{3}\in \mathcal{GP}(\xi)$. Then we have the following two commutative diagrams by base change:
\begin{center}$\xymatrix{&G_{3}\ar@{=}[r]\ar[d]&G_{3}\ar[d]\\
\Sigma^{-1}A\ar[r]\ar@{=}[d]&W\ar[r]\ar[d]&Q\ar[r]\ar[d]&
A\ar@{=}[d]\\
\Sigma^{-1}A\ar[r]&K\ar[r]\ar[d]&G_{0}\ar[r]\ar[d]&A\\
&\Sigma G_{3}\ar@{=}[r]&\Sigma G_{3}\\
}$
\end{center}

\begin{center}$\xymatrix{&G_{3}\ar@{=}[r]\ar[d]&
G_{3}\ar[d]\\
B\ar[r]\ar@{=}[d]&G_{1}'\ar[r]\ar[d]&W\ar[r]\ar[d]&
\Sigma B\ar@{=}[d]\\
B\ar[r]&G_{1}\ar[r]\ar[d]&K\ar[r]\ar[d]&\Sigma B\\
&\Sigma G_{3}\ar@{=}[r]&\Sigma G_{3}\\
}$
\end{center}
Since $\xi$ is closed under base change, we get that the triangles $G_{3}\rightarrow W\rightarrow K\rightarrow \Sigma G_{3}$ and $B\rightarrow G_{1}'\rightarrow W\rightarrow \Sigma B$ are in $\xi$.
Applying the functor $\Hom_{\mathcal{C}}(\mathcal{P}(\xi),-)$ to the above two diagrams, by snake lemma and Lemma \ref{lem1} we have that the triangles $W\rightarrow Q\rightarrow A\rightarrow \Sigma W$ and $G_{3}\rightarrow G_{1}'\rightarrow G_{1}\rightarrow \Sigma G_{3}$ are in $ \xi$. Because both $G_{1}$ and $G_{3}$ are in $\mathcal{GP}(\xi)$, $G_{1}'$ is also in $\mathcal{GP}(\xi)$ by \cite[Theorem 3.11]{AS04}. Connecting the triangles $B\rightarrow G_{1}'\rightarrow W\rightarrow \Sigma B$ and $W\rightarrow Q\rightarrow A\rightarrow \Sigma W$, we get the second desired $\xi$-exact sequence. \hfill$\Box$

In particular, we have the following corollary.

\begin{corollary}\label{cor1} Let $ G_{1}\rightarrow G_{0}\rightarrow A\rightarrow \Sigma G_{1}$ be in $\xi$ with $G_{1},G_{0} \in \mathcal{GP}(\xi)$. Then there exist the following triangles:
$$ P\rightarrow G_{0}'\rightarrow A\rightarrow \Sigma P$$
and $$G_{1}'\rightarrow Q\rightarrow A\rightarrow \Sigma G_{1}'$$
 in $\xi$ where $P,Q\in \mathcal{P}(\xi)$ and $G_{0}',G_{1}'\in \mathcal{GP}(\xi)$.
\end{corollary}

\begin{proposition} Let $\mathcal{C}$ be a triangulated category with enough $\xi$-projectives, $A\in \mathcal{C}$,  and $n$ be a non-negative integer. Then the following statements are equivalent:

(1) $\xi$-$\Gpd(A)\leq n$.

(2) For every $0\leq i\leq n$, there is an $\xi$-exact sequence $$0\rightarrow P_{n}\rightarrow \cdots \rightarrow P_{i+1}\rightarrow G\rightarrow P_{i-1}\rightarrow \cdots\rightarrow P_{0}\rightarrow A\rightarrow 0$$
 with $P_{j}\in \mathcal{P}(\xi)$ for all $0\leq j\leq n$, $j\neq i$, and $G\in \mathcal{GP}(\xi)$.

(3) For every $0\leq i\leq n$, there is an $\xi$-exact sequence $$0\rightarrow G_{n}\rightarrow \cdots \rightarrow G_{i+1}\rightarrow P\rightarrow G_{i-1}\rightarrow \cdots\rightarrow G_{0}\rightarrow A\rightarrow 0$$
 with $G_{j}\in \mathcal{GP}(\xi)$ for all $0\leq j\leq n$, $j\neq i$, and $P\in \mathcal{P}(\xi)$.
\end{proposition}

\noindent{\bf Proof} The case $n=0$ is trivial. We may assume $n\geq 1$.

(1)$\Rightarrow $(2) we proceed by induction on $n$. Suppose $\xi$-$\Gpd(A)\leq 1$. Then there exists a triangle $G_{1}\rightarrow G_{0}\rightarrow A\rightarrow \Sigma G_{1}$ in $\xi$ with $G_{0},G_{1}\in \mathcal{GP}(\xi)$. By Corollary \ref{cor1}, we get the triangles $P\rightarrow G_{0}'\rightarrow A\rightarrow \Sigma P$ and $G_{1}'\rightarrow Q\rightarrow A\rightarrow \Sigma G_{1}'$ in $\xi$ with $P,Q\in \mathcal{GP}(\xi)$ and $G_{0}',G_{1}'\in \mathcal{GP}(\xi)$.

Now suppose $n\geq 2$. Then there exists an $\xi$-exact sequence
$$0\rightarrow G_{n}\rightarrow  G_{n-1}\rightarrow \cdots\rightarrow  G_{1}\rightarrow G_{0}\rightarrow A\rightarrow 0$$
with $G_{i}\in \mathcal{GP}(\xi)$ for all $0\leq i\leq n$. Applying Proposition \ref{prop1} to the relevant $\xi$-exact sequence $0\rightarrow K\rightarrow G_{1}\rightarrow G_{0}\rightarrow A\rightarrow 0$, we get an $\xi$-exact sequence $0\rightarrow K\rightarrow G_{1}'\rightarrow P_{0}\rightarrow A\rightarrow 0$ with $G_{1}'\in \mathcal{GP}(\xi)$ and $P_{0}\in \mathcal{P}(\xi)$, which yields an $\xi$-exact sequence $$0\rightarrow G_{n}\rightarrow  G_{n-1}\rightarrow \cdots\rightarrow  G_{2}\rightarrow  G_{1}'\rightarrow P_{0}\rightarrow A\rightarrow 0.$$
Taking into account the relevant $\xi$-exact sequence $$0\rightarrow G_{n}\rightarrow  G_{n-1}\rightarrow \cdots\rightarrow  G_{2}\rightarrow  G_{1}'\rightarrow L\rightarrow 0,$$ it follows that $\xi$-$\Gpd(L)\leq n-1$. By the induction hypothesis, there exists an $\xi$-exact sequence $$0\rightarrow P_{n}\rightarrow \cdots \rightarrow P_{i+1}\rightarrow G\rightarrow P_{i-1}\rightarrow \cdots\rightarrow P_{1}\rightarrow L\rightarrow 0$$
 with $P_{j}\in \mathcal{P}(\xi)$ for all $1\leq j\leq n$, $j\neq i$, and $G\in \mathcal{GP}(\xi)$.   Now one can paste the above $\xi$-exact sequence and the triangle $L\rightarrow P_{0}\rightarrow A\rightarrow \Sigma L$ together to obtain the desired $\xi$-exact sequence.

 (2)$\Rightarrow $(1) and (3)$\Rightarrow $(1) are clear.

 (1)$\Rightarrow $(3) Suppose $\xi$-$\Gpd(A)\leq n$. Then there exists an $\xi$-exact sequence
$$0\rightarrow G_{n}\rightarrow  G_{n-1}\rightarrow \cdots\rightarrow  G_{1}\rightarrow G_{0}\rightarrow A\rightarrow 0$$
with $G_{i}\in \mathcal{GP}(\xi)$ for all $0\leq i\leq n$. For every $0\leq i<n$, considering the relevant $\xi$-exact sequence
$$0\rightarrow G_{n}\rightarrow  G_{n-1}\rightarrow \cdots\rightarrow  G_{i+1}\rightarrow G_{i}\rightarrow M\rightarrow 0,$$ it follows that $\xi$-$\Gpd(M)\leq n-i$. By the proof of (1)$\Rightarrow $(2), we get an $\xi$-exact sequence
$$0\rightarrow G_{n}'\rightarrow  G_{n-1}'\rightarrow \cdots\rightarrow  G_{i+1}'\rightarrow P\rightarrow M\rightarrow 0$$ with $G_{i}'\in \mathcal{GP}(\xi)$ and $P\in \mathcal{P}(\xi)$. So we obtain the $\xi$-exact sequence $$0\rightarrow G_{n}'\rightarrow \cdots \rightarrow G_{i+1}'\rightarrow P\rightarrow G_{i-1}\rightarrow \cdots\rightarrow G_{0}\rightarrow A\rightarrow 0.$$

Now we only need to prove the result for $i=n$. Applying Corollary \ref{cor1} to the relevant triangle $G_{n}\rightarrow  G_{n-1}\rightarrow L\rightarrow  \Sigma G_{n}$, we get the triangle $P\rightarrow  G_{n-1}'\rightarrow L\rightarrow  \Sigma P$ with $G_{n-1}'\in \mathcal{GP}(\xi)$ and $P\in \mathcal{P}(\xi)$. Thus we obtain the desired $\xi$-exact sequence $$0\rightarrow P \rightarrow G_{n-1}'\rightarrow G_{n-2}\rightarrow \cdots\rightarrow G_{0}\rightarrow A\rightarrow 0.$$ \hfill$\Box$

\end{document}